\documentclass[12pt]{article}
\usepackage{amsmath}
\usepackage{amssymb}
\usepackage[english]{babel}

\newtheorem{Theorem}{Theorem}[section]
\newtheorem{Lemma}[Theorem]{Lemma}
\newtheorem{Corollary}[Theorem]{Corollary}
\newtheorem{Proposition}[Theorem]{Proposition}
\newtheorem{Note}[Theorem]{Remark}
\newtheorem{Example}[Theorem]{Example}

\def\N{\hbox{${\mathbb N}$}}

\def\ep{\hfill{\vbox to 7pt{\hbox to 7pt{\vrule height 7pt width 7pt}}}}
\addtocounter{section}{-1}

\title{On Bernstein algebras satisfying chain conditions II}

\author{Fouad Zitan}

\begin{document}
	\maketitle 
\centerline{Universit\'e Abdelmalek Essa\^adi, Facult\'e des Sciences} 
\centerline{D\'epartement de Math\'ematiques, B.P.
2121, T\'etouan, Morocco}
\centerline{e-mail: zitanf@yahoo.fr}
	\vspace*{0,7cm}
	\centerline{\it To the memory of my wonderful teacher, Professor Artibano Micali $(1931$-$2011)$}
 
\baselineskip0.6cm
\vspace*{0,4cm}\vspace*{0,2cm}
\begin{abstract} 
	\vspace*{0,1cm}
Following a previous work with  Boudi, we continue to  investigate Bernstein algebras
satisfying chain conditions. First, it is shown that  a Bernstein algebra $(A, \omega)$ with ascending  or   descending chain condition on subalgebras
is finite-dimensional. We also prove  that  $A$ is N\oe
therian (Artinian) if and only if
 its barideal $N=\ker(\omega)$ is. Next,  as a generalization of Jordan and nuclear Bernstein algebras, we study whether a N\oe therian (Artinian) Bernstein algebra $A$ with a locally nilpotent
 barideal $N$ is finite-dimensional. The response is affirmative in the N\oe therian case, unlike in the Artinian case. This question is closely related to a result by Zhevlakov on general locally nilpotent nonassociative algebras that are N\oe therian, for which we give a new proof. In particular, we derive that a commutative nilalgebra of nilindex 3 which is N\oe therian or Artinian is finite-dimensional. Finally, we improve and extend some results of Micali and Ouattara to the N\oe therian and Artinian cases. 
\end{abstract}

\vspace*{0,3cm}
{\small {\it 2010 Mathematics Subject Classification}:
Primary 17D92; Secondary  17C10, 17A01.\\

\hspace*{2,4cm} N\oe therian, nuclear algebra.}

\section{\bf Introduction}

\hspace*{0,5cm} Bernstein algebras form a class of
nonassociative algebras whose origin lies in genetics. Historically, they have been introduced by Lyubich \cite{Lyu} and Holgate \cite{Ho} as an algebraic formulation of the problem of classifying the stationary evolution operators in genetics. In this way, Bernstein algebras represent populations reaching the equilibrium
after the first generation. Since then the theory has evolved into an independent branch of nonassociative algebras,  and many researches have been done on the topic from various points of view (see, for instance, \cite{Be, Co, Gonza1, Ly, Mi, Ouattara, Wa, Wo}).\\

 One of the main questions on the structure of Bernstein
algebras, posed by  Lyubich and solved by  Odoni and
 Stratton \cite{Od} as well as by  Baeza \cite{Ba} and Grishkov \cite{Gr}, says that the barideal of a
finite-dimensional nuclear Bernstein algebra is nilpotent. The
analogous question in the finitely generated case was proposed by
 Grishkov \cite{Gr}, and its affirmative solution was settled
by Peresi \cite{Pe} and Krapivin \cite{Kra}. Consequently, finitely generated nuclear Bernstein algebras must be finite-dimensional,
as was directly established by Suazo \cite{Su} using a different approach.\\

On the other hand, it is well known that one of the most satisfactory developments of the theory of some varieties of algebras, as associative, Jordan, alternative and Lie algebras, is the structure theory of algebras with chain conditions, and there is presently a substantial bibliography in this subject. Concerning  Bernstein algebras, a detailed treatment
 was done in \cite{BZ} for Bernstein algebras satisfying
chain conditions on ideals. Among many other results in that paper, it was especially proved that for a Bernstein algebra which is Jordan or nuclear, each of the N\oe therian and Artinian hypotheses implies finite-dimensionality of the algebra and so nilpotency of the barideal. Hence, the Lyubich conjecture is still valid in the N\oe therian and
Artinian cases.\\ 

In this present article we pursue the study of Bernstein algebras satisfying conditions initiated  in \cite{BZ}. After a first
section devoted to preliminaries, we show in Section 2 that  an arbitrary Bernstein algebra $A$  satisfying the ascending  or descending chain condition subalgebras
is finite-dimensional.  Then we prove in Section 3  that a Bernstein algebra $A$ is N\oe therian (Artinian) if and only if
its barideal $N=\ker(\omega)$ is,  thus generalizing an early result by Krapivin \cite{Kra} about the finitely-generated case. Section 4  deals with Bernstein algebras having locally nilpotent
barideals, as an extension of both Jordan and nuclear Bernstein algebras. Specifically, we study whether a N\oe therian (Artinian) Bernstein algebra $A$ with a locally nilpotent
barideal $N$ is finite-dimensional. The answer is positive in the N\oe therian case and negative in the Artinian case. This question is connected with a result due to Zhevlakov \cite{Zh} on general locally nilpotent nonassociative algebras for which we provide an independent proof. As a special case, we derive that a commutative nilalgebra of nilindex 3 which is N\oe therian or Artinian is finite-dimensional.
In the final section,  we intend to ameliorate and extend  some results of Micali and Ouattara  \cite{Mi} to the N\oe therian and Artinian cases.\\ 
Various examples are presented along this work to serve as motivation and illustration for our results.

\section{\bf Preliminaries}
\quad In this section we briefly summarize notation, terminology and
classical properties about Bernstein algebras and arbitrary nonassociative algebras.
 Throughout this paper, we will fix an infinite ground field $K$
of characteristic different from 2 and 3, and let $A$ be an algebra over $K$, not
necessarily associative or finite-dimensional. If there exists a nonzero  homomorphism of algebras
$\omega:A \rightarrow K$, then the ordered pair $(A, \omega)$
is called a {\it baric algebra} and $\omega$ is its {\it weight
function}. For every $e\in A$ with $\omega(e)\neq 0$, we have
$A=Ke\oplus N$, where $N=\ker(\omega)$ is an ideal of $A$,
called the {\it barideal} of $A$. A {\it baric ideal} of $A$ 
is an ideal $I$ of $A$ with $I\subseteq N$. Then the quotient algebra $A/I$ is a baric algebra with weight function $\overline{\omega}$  defined by $\overline{\omega}(x+I)=\omega(x)$.\\
\hspace*{0,5cm} A {\it Bernstein algebra} is a commutative baric algebra $(A, \omega)$
satisfying the identity $(x^2)^2=\omega(x)^2x^2$. A Bernstein algebra has a unique weight function $\omega$. If $x\in A$ and
$\omega(x)=1$, then $e=x^2$ is a nontrivial idempotent of $A$ which gives rise to the {\it Peirce decomposition} $A=Ke\oplus U\oplus V$, where
$N=U\oplus V$, and
\begin{equation}
U=\{u\in A\; /\; eu=\frac 12\ u\},\quad  V=\{v\in A\; /\;
ev=0\}.
\end{equation}
Besides, the Peirce components multiply according to  the relations
\begin{equation}
U^2\subseteq V, \; UV\subseteq U, \; V^2\subseteq U, \; UV^2=0.
\end{equation}
A Bernstein algebra $A=Ke\oplus U\oplus V$ is never unital unless if $\dim(A)=1$, and cannot be associative except when $U=0$. However, Bernstein algebras may be power-associative, that is, if each element generates  an associative subalgebra. Recall that a commutative
algebra $A$ is a {\it Jordan algebra} if the identity $x(x^2
y)=x^2(xy)$ holds in $A$. Bernstein-Jordan algebras play a crucial role in the theory of Bernstein algebras. It is  well-known  that the following four conditions are equivalent for a Bernstein
algebra $A=Ke\oplus U\oplus V$ (see, for instance, 
\cite{Gonza1,  Wa}):\\
(a) $A$ is a Jordan algebra. \qquad\qquad\quad \ (b) $A$ is power-associative.\\
(c) $x^3=\omega(x)x^2$ for all $x\in A$. \qquad \qquad (d) $V^2=0$ and $(uv)v=0$ for all $u\in U$ and $v\in V$. \\ 

Therefore,  the
elements of the barideal $N=\ker(\omega)$ in a Bernstein-Jordan algebra $(A, \omega)$ satisfy $x^3=0$ and so the Jacobi
identity
\begin{equation}
(xy)z+(yz)x+(zx)y=0.
\end{equation} 
\hspace*{0,5cm} An important tool in Bernstein algebras is the ideal $ann_U(U)=\{u\in U \; /\; uU=0\}$ of
$A$ wich is independent of the selected idempotent $e$ and satisfies
$ann_U(U)(U\oplus U^2)=0$ and  $V^2\subseteq ann_U(U)$ (see, for instance,
\cite[Theorem 3.4.19]{Ly}). It should be remarked the
fundamental position played by this ideal $ann_U(U)$ in the connection between
Bernstein algebras and Jordan algebras, since the quotient algebra
$A/ann_U(U)$ is a Bernstein-Jordan algebra (see, for example, \cite{Gonza1, Hen, Mi}). A Bernstein algebra $A$ is called {\it nuclear} if
$A^2=A$, or equivalently, $U^2=V$ for an arbitrary idempotent $e$ ; in this case, we have
$ann_U(U)N=0$. Every Bernstein
algebra $A=Ke\oplus U\oplus V$ gives rise to  a nuclear Bernstein subalgebra $A^2=Ke\oplus U\oplus U^2$.  Further information about algebraic properties of Bernstein algebras, as well as their possible genetic interpretation can be found in \cite{Ly, Ouattara, Reed,  Wo}.\\
 \hspace*{0,5cm} We let now $A$ be an arbitrary algebra over $K$. Following
 the notation of \cite{russe}, we will consider the {\it powers}
 $A^i$ and the {\it right principal powers} $A^{<i>}$ of $A$
 defined recursively by $A^1=A^{<1>}=A, \; A^i=\sum\limits_{r+s=i}
  A^r A^s$ and $A^{<i>}=A^{<i-1>}A$. The algebra $A$ is called
  {\it nilpotent} if $A^n=0$ for some $n$, and {\it right
  nilpotent} if $A^{<n>}=0$. It is well known that these two
  notions of nilpotency are equivalent in the commutative case. We
  define the plenary powers $A^{(i)}$ of $A$ by setting
  $A^{(1)}=A^2$ and $A^{(i)}=(A^{(i-1)})^2$. The algebra $A$ is
  said to be {\it solvable} when $A^{(n)}=0$ for some $n$.\\
  \hspace*{0,5cm} On the other hand, $A$ being a commutative algebra, the
 linear mappings $L_a: A \rightarrow A$ defined by $L_a(x)=ax$  generate a
 subalgebra of $End_K(A)$, denoted by ${\mathcal M}_*(A)$ and
 called the {\it multiplication ideal} of $A$. The subalgebra of
 $End_K(A)$ generated by ${\mathcal M}_*(A)$ and the identity
 endomorphism id$_A$, will be denoted by ${\mathcal M}(A)$ and will
 be called the {\it multiplication algebra} of $A$. If $B$ is a
 subalgebra of $A$, we write ${\mathcal M}_*^A(B)$ for the
 subalgebra of ${\mathcal M}_*(A)$ generated by the operators
 $L_b$, where $b\in B$. The unital
 algebra ${\mathcal M}^A(B)$ is defined analogously.\\
 \hspace*{0,5cm} For any subset $S\subseteq A$, we will adopt the notations $<S>$
 and $K\langle S\rangle$, which mean respectively the subspace of
 $A$ spanned by $S$ and the free unital associative
 (noncommutative) algebra over $K$ generated by $S$. 
 Since the ideal of $A$ generated by $S$ consists of finite sums of elements $f(x)$,
 where $f\in {\mathcal M}(A)$ and
 $x\in S$, it is customary to denote it by ${\mathcal M}(A)S$. \\
  \hspace*{0,5cm} Returning to Bernstein algebras, recall that in a Bernstein
 algebra, the  principal powers $N^{<i>}$ are ideals \cite[page 113]{Ly}.  Moreover, The barideal
  $N$ satisfies the equation
 $(x^2)^2=0$, but is not in general nilpotent. However, $N$ is
 always solvable, since $N^{(3)}=0\;$ \cite[Theorem 2.11]{Be} (see, also, \cite{Jac}).

 \section{\large Chain conditions for subalgebras}

\hspace*{0,5cm} In our preceding work \cite{BZ} it has been proved that a Bernstein algebra $A$ that
is Jordan or nuclear is necessarily finite-dimensional whenever it
is N\oe therian or Artinian. In addition, a counter-example  was given to show that the hypothesis that $A$ be Jordan or nuclear  is essential in this result.
In the following we are going to relax the Jordan and nuclear assumptions  in order to state a result for general Bernstein algebras
satisfying the ascending (descending) condition for subalgebras
instead of ideals. 
An arbitrary algebra $A$  satisfies the ascending chain condition a.c.c.
(descending chain condition d.c.c.) on subalgebras if it has no infinite strictly ascending (descending) chains of subalgebras. It is easily seen that the a.c.c. (d.c.c.) for subalgebras is equivalent to the maximal (minimal) condition for
subalgebras, that is, every non-empty set of subalgebras has a
maximal (minimal) element. Moreover, in an algebra satisfying
the a.c.c. for subalgebras, all subalgebras are finitely generated.
In the literature of general nonassociative algebras, there are some
results treating the maximal condition for subalgebras. For
instance, Kubo constructed in \cite{Kubo} infinite-dimensional
associative, Jordan and Lie algebras satisfying
the maximal condition for subalgebras (see also \cite{Amayo}). In \cite[Theorem 3, page
91]{russe} it is established that a Jordan nil-algebra satisfying
the maximal condition for subalgebras is nilpotent, and therefore
finite-dimensional. For Bernstein algebras, we may formulate the
following result which is valid for both the a.c.c and d.c.c.
conditions on subalgebras.

\begin{Theorem} For a Bernstein algebra $A$, the following conditions are equivalent:\\
{\rm (i)} $A$ satisfies the a.c.c. (d.c.c.) condition for  subalgebras; \\
{\rm (ii)} $A$ satisfies the a.c.c. (d.c.c.) condition for  subalgebras contained in $N=\ker(\omega)$;\\
{\rm (iii)} $A$ is finite-dimensional.
\end{Theorem}

{\it Proof. } It is enough to demonstrate that (ii) implies (iii). If (ii) holds, then $A$ satisfies a fortiori the a.c.c.
(d.c.c.) condition for ideals contained in $\ker(\omega)$, and hence it is
N\oe therian (Artinian) in view of \cite[Proposition 2.1]{BZ}. It
follows from \cite[Proposition 3.1]{BZ} that the Bernstein-Jordan
algebra $A/ann_U(U)$ is finite-dimensional. Now, since
$\left(ann_U(U)\right)^2=0$, every subspace of $ann_U(U)$ is a
subalgebra of $A$ contained in $\ker(\omega)$. It follows from the
hypothesis that $ann_U(U)$ is finite-dimensional, which completes
the proof. \ep

\section{\large The barideal of a Bernstein algebra}

\hspace*{0,5cm} Krapivin established in
\cite{Kra} that a Bernstein algebra $(A, \omega)$ is finitely
generated if and only if its barideal $N=\ker(\omega)$ is finitely
generated (as an algebra). Hence, it is legitimate to ask the analogous question for the N\oe therian and Artinian
cases. An arbitrary algebra is said to be N\oe therian
(Artinian) if it satisfies the ascending chain condition a.c.c.
(descending chain condition d.c.c.) on ideals, that is, every
ascending (descending) sequence of ideals is stationary. Before embarking in this direction, we require some
preparation. The key ingredient is the deep   link  exhibited in \cite{BZ}   
between Bernstein algebras and modules over
associative (noncommutative) algebras. In details, let $A=Ke\oplus U\oplus V$ be a Bernstein algebra, and
consider the free unital associative (noncommutative) algebra
$K\langle V\rangle$ generated by the set $V$. Then  the ideal
$ann_U(U)$ becomes a left module over $K\langle V\rangle$ by
setting
$$(v_1*\ldots *v_k).u=v_1(\ldots (v_ku)\ldots ), \mbox{ for all }
v_1, \ldots, v_k\in V \mbox{ and } u\in ann_U(U).$$ This $K\langle
V\rangle$-module $ann_U(U)$ contains many information on the
Bernstein algebra $A$. For instance, the submodules of this module
are just the ideals of $A$ contained in $ann_U(U)$. Moreover, the
finiteness behavior of the Bernstein algebra $A$ was studied with
much benefit in terms of its attached $K\langle V\rangle$-module
$ann_U(U)$. In particular, it was established that the Bernstein
algebra $A$ is finitely generated (resp. N\oe therian, Artinian)
if and only if $A/ann_U(U)$ is finite-dimensional and the
$K\langle V\rangle$-module $ann_U(U)$ is finitely generated (resp.
N\oe therian, Artinian).

\hspace*{0,5cm} Now, we are in a position to prove the following result
which extends the Krapivin theorem \cite{Kra} to both the N\oe therian and Artinian contexts. 

\begin{Theorem} Let $A=Ke\oplus U\oplus V$ be a Bernstein algebra
	with barideal $N=\ker(\omega)=U\oplus V$. Then the following
	conditions are equivalent: \\
	{\rm (i)}  $A$ is N\oe therian (Artinian); \\
	{\rm (ii)}  $N$ is N\oe therian (Artinian).
\end{Theorem}

{\it Proof. } The implication ${\rm (ii)} \Rightarrow {\rm (i)}$ is trivial, since a
	Bernstein algebra $A$ is N\oe therian (Artinian) if and only if it
	satisfies a.c.c. (d.c.c.) on baric ideals of $A$
	\cite[Proposition 2.1]{BZ}.\\
	${\rm (i)} \Rightarrow {\rm (ii)}$: Assume that $A$ is N\oe therian. Then by
	\cite[Proposition 3.4]{BZ}, $A/ann_U(U)$ is finite-dimensional and
	the $K\langle V\rangle$-module $ann_U(U)$ is N\oe therian. If $I$
	is an ideal of $N$, then $I\cap ann_U(U)$ is an ideal of $A$
	contained in $ann_U(U)$, because $ex=\frac 12 x$ for all $x\in
	I\cap ann_U(U)\subseteq U$. Hence, $I\cap ann_U(U)$ is a submodule
	of the
	$K\langle V\rangle$-module $ann_U(U)$.\\
	Now, let $(I_n)$ be an increasing sequence of ideals of $N$. Then
	$(I_n\cap ann_U(U))$ is an increasing sequence of submodules of
	the $K\langle V\rangle$-module $ann_U(U)$, and
	$(I_n+ann_U(U))/ann_U(U)$ is an increasing sequence of subspaces
	of the quotient space $A/ann_U(U)$. Then both chains must
	stabilize, and by a standard argument, the sequence $(I_n)$ is
	stationary. \\The Artinian case is treated analogously. \ep

\section{\large Bernstein algebras and locally nilpotent nonassociative algebras}

\hspace*{0,5cm} Jordan and nuclear Bernstein algebras are important types of Bernstein algebras. The barideal $N=\ker(\omega)$ in a Bernstein-Jordan algebra $(A, \omega)$ satisfies the identity $x^3=0$, hence by \cite[page 114]{russe} $N$ is {\it locally nilpotent},
that is, every finitely generated subalgebra of $N$ is
nilpotent (see, also, \cite{Co}). Let us explain the similar fact for nuclear Bernstein algebras:
 
\begin{Proposition} Let $A$ be a nuclear Bernstein algebra. Then
the barideal $N$ of $A$ is locally nilpotent.
\end{Proposition}

{\it Proof. } Consider the Bernstein-Jordan algebra
$\overline{A}=A/ann_U(U)$ and let $\pi : A
\longrightarrow A/ann_U(U)$ be the canonical surjection. We know that the barideal $N/ann_U(U)$
of the Bernstein-Jordan algebra $A/ann_U(U)$ is locally nilpotent.
Let the subalgebra $S$ of $A$ generated by the elements $a_1,
\dots, a_n$. Then the subalgebra $T=\pi(S)$ of $A/ann_U(U)$
generated by $\pi(a_1), \dots, \pi(a_n)$ is nilpotent,
say $T^{<k>}=0$. It follows that
 $S^{<k>}\subseteq ann_U(U)$. Hence,
$S^{<k+1>}=S^{<k>}S\subseteq ann_U(U)S\subseteq ann_U(U)N=0$, which
means that $S$ is nilpotent. It follows that  $N$ is locally nilpotent. \ep

\begin{Note} {\rm Let $A=Ke\oplus U\oplus V$ be an arbitrary Bernstein
algebra. Then the subspace $U\oplus U^2$ is an ideal of $A$ which
is locally nilpotent. To be convinced, it suffices to consider the
nuclear Bernstein subalgebra $A^2=Ke\oplus U\oplus U^2$ whose
barideal is $U\oplus U^2$.}
\end{Note}

In virtue of \cite[Theorem 2.3]{BZ}, for a Bernstein algebra which is Jordan or
nuclear, each of the N\oe therian and Artinian conditions implies
finite-dimensionality of the algebra and so nilpotency of the
barideal. This result suggests us to
  raise the following more general question: Let $(A, \omega)$ be a
  Bernstein algebra such that its barideal $N=\ker(\omega)$ is
  locally nilpotent. If $A$ is N\oe therian or
  Artinian, is it finite-dimensional?   In the light of our  Theorem 3.1,
  the above question has a closed link with the
  following question on general locally nilpotent nonassociative
  algebras, which seems to be of independent
interest: Is a locally nilpotent algebra which is  N\oe therian or Artinian finite-dimensional? Searching in the wide literature
  of general nonassociative algebras, we found a noteworthy
  article \cite{Zh} published in 1972 by the eminent algebraist Zhevlakov (1939-1972) after his death, 
  which gives a positive answer to the N\oe therian case and
  constructs a counter-example to the Artinian case. For the sake of completeness, we  provide below  an alternative 
  proof of this result which is substantially different from the proof of Zhevlakov mentioned in \cite[Note 1]{Zh}, and that we have done before discovering Zhevlakov's paper.

\begin{Theorem} Let $N$ be a nonassociative (possibly noncommutative)
algebra over a field $K$. Assume that $N$ satisfies the ascending
chain condition on two-sided ideals. If $N$ is locally nilpotent,
then $N$ is finite-dimensional.
\end{Theorem}

{\it Proof. } Let the ideal $N$ be generated by the elements $e_1,\dots,e_r$:
$N={\cal M}(N)e_1+\dots+{\cal M}(N)e_r$. We denote by $F$ the
subalgebra generated by $e_1,\dots,e_r$. Then $F$ is nilpotent, and
by a straightforward argument, $F$ is
finite-dimensional.  Choose
$m\geq 2$ such that the power $F^m=0$. Clearly, $N=F+N^2$, and by a simple induction one may show that 
\begin{equation} N^i\subseteq F^i+N^{i+1}\mbox{ for each } i\geq 1\end{equation}
 Indeed, if the inclusion (4) is true for each $j\leq i$, let us show that $N^{i+1}\subseteq F^{i+1}+N^{i+2}$. We have: 
$$N^{i+1}=\sum\limits_{j_1+j_2=i+1}
N^{j_1}N^{j_2}\subseteq \sum\limits_{j_1+j_2=i+1} (F^{j_1}+N^{j_1+1})(F^{j_2}+N^{j_2+1}).$$
Now, the following relations are easy to verify: $$F^{j_1}F^{j_2}\subseteq F^{j_1+j_2}=F^{i+1}, \;\;\;\;\;
F^{j_1}N^{j_2+1}\subseteq N^{j_1}N^{j_2+1}\subseteq N^{j_1+j_2+1}=N^{i+2},$$
$$N^{j_1+1}F^{j_2}\subseteq N^{j_1+1}N^{j_2}\subseteq N^{j_1+1+j_2}=N^{i+2},\;\;\;\;\; N^{j_1+1}N^{j_2+1}\subseteq N^{j_1+j_2+2}=N^{i+3}\subseteq N^{i+2},$$
from which we get the desired inclusion $N^{i+1}\subseteq F^{i+1}+N^{i+2}$.\\
As a consequence of (4), it follows
 that $N^i\subseteq F+N^{i+1}$, for each $i\geq 1$. \\Hence,
$N=F+N^2=F+N^3=\dots=F+N^m$. Now, consider $a_1,\dots,a_t$ in $N^m$ such
that, $N^m={\cal M}(N)a_1+\dots+{\cal M}(N)a_t$. Without loss of
generality, we can assume that each $a_j$ is  a nonassociative
product $x_1\dots x_m$ (with some distribution of parentheses) of $m$
factors $x_i\in N$. Writing $x_k=u_k+v_k$ $(k=1,\dots,m)$, where
$u_k\in F$, $v_k\in N^m$ and using the fact that $F^m=0$, we obtain
that each $a_j$ is a sum of elements of the form $y_1\dots y_m$ (with
some distribution of parentheses), where $y_{1},\dots,y_{m}\in N$
such that at least one of them, say $y_{d}$, belongs to $N^m$.
Decomposing $y_{d}$ into $N^m={\cal M}(N)a_1+\dots+{\cal M}(N)a_t$,
we get that $y_{1}\dots y_{d}\dots y_{m}$ is a sum of products
$h_1\dots h_{s-1}a_ih_{s+1}\dots h_p$ (with some distribution of
parentheses), where $h_1,\dots,h_{s-1}, h_{s+1},\dots, h_p\in N$ and
$i\in \{1,\dots,t\}$. As a consequence, each $a_j$ can be
expressible in the form
\begin{equation}
a_j =\sum h_1^{j,i}\dots h_{s_{j,i}-1}^{j,i}\;a_i\;
h_{s_{j,i}+1}^{j,i}\dots h_{p_{j,i}}^{j,i},
\end{equation}
where $h_k^{j,i}\in N,\; p_{j,i}\geq s_{j,i}\geq 1.$ \\
The subalgebra $H$ generated by the finite set $\{a_1,\dots, a_t\} \cup
\{h_k^{j,i}\}$ is nilpotent. Let $H^l=0$ ($l\geq 2$). We may now
replace $a_i$ by its expression ($l$ times) in (5) to conclude
that $a_j=0$. This implies that $N^m=0$, and therefore $N=F$, which completes the proof. \ep\\

A special case of locally nilpotent algebras are commutative
nilalgebras of nilindex at most 3, whose natural  examples are
barideals of Bernstein-Jordan algebras and train algebras of rank 3 (see, for instance, \cite{Zitan}). They satisfy the identity
$x^3=0$ and so also the Jacobi identity $(xy)z+(yz)x+(zx)y=0$. It
is well known that they are Jordan algebras (see \cite[Lemma
1]{Ba}, \cite[Lemma 2.2]{Jacobi}, \cite{Gut}, \cite[page 114]{russe}). They appear in the litterature as {\it Jacobi-Jordan algebras} in
\cite{Agore,Jacobi}  and as {\it mock-Lie algebras} in \cite{Pasha}. In addition, any such an algebra $N$ is solvable
and $N^{(4)}=0$ \cite[Lemma 3.1]{Zitan}. On the other hand, since such algebras are locally nilpotent, we infer the following
immediate consequence of Theorem 4.3.

\begin{Corollary} Let $N$ be a commutative algebra satisfying the identity $x^3=0$. If $N$ satisfies the ascending
chain condition on ideals, then $N$ is finite-dimensional.
\end{Corollary}

The above  corollary is certainly not new, and it is quite easy to prove it directly. Indeed, since the ideal $N$ is a N\oe therian solvable  Jordan
algebra, it follows from a result of Medvedev and Zelmanov \cite{Med} that $N$ is nilpotent. Hence, by a simple reasoning, one may prove that $N$ is finite-dimensional.\\

\hspace*{0,5cm} The analog of Theorem 4.3 for the descending chain condition
is false, as already shown by Zhevlakov \cite{Zh}: 

\begin{Example} {\rm (Zhevlakov) Let $N$ be a countably-dimensional  algebra with basis
		$\left\{e_n\right\}_{n\in \N^*}$ and nonzero products
		 $$ e_ie_j=e_{\min(i,j)-1} \;\mbox{ for }i, j\geq 2.$$  By \cite[Note 2]{Zh}, $N$ is a  locally nilpotent commutative algebra which is  Artinian and satisfying $N^2=N$.}
\end{Example}

We offer  another counter-example below: 

\begin{Example} {\rm Let $N$ be an infinite dimensional commutative algebra with basis
$\left\{e_n\right\}_{n\in \N^*}$ and nonzero multiplication table
given by $ e_n^2=e_{n-1}$ for $n\geq 2$.  \\
Let $S$ be a nonzero subalgebra of $N$, and let $x=\alpha _1e_1+\dots+\alpha _ke_k\in S$, with
$\alpha_k\neq 0$. By considering the plenary powers $x^{[r]}\;
(1\leq r\leq k)$  defined by $x^{[1]}=x$ and
$x^{[r]}=(x^{[r-1]})^2$, a simple calculation gives $x^{[2]}=\alpha _2e_1+\dots+\alpha
_ke_{k-1},\;\ldots, x^{[k-1]}=\alpha_{k-1}e_1+\alpha_k e_2,
\;x^{[k]}=\alpha_{k}e_1$. It follows that $e_1,\dots, e_k\in S$. If $n=\max\,\{k \; / \; x\neq 0 \} $ is finite, then
$S=<e_1, \ldots, e_n>$, which is also an ideal of $N$. In the opposite case, we get $S=N$. Therefore, each proper subalgebra (ideal) of $N$ coincides with a subspace $<e_1, \ldots, e_n>$ for some $n\geq 1$. Clearly, no infinite sequence of ideals of $N$ can exist, and so $N$ is Artinian. On the other hand, it is not hard to observe that $N$
is locally nilpotent, since every finitely generated subalgebra $<e_1, \ldots, e_n>$ is nilpotent (of index $n+1$). Obviously, $N$ is not N\oe therian because the
ascending chain
$(<e_1, \ldots, e_n>)_{n\geq 1}$
of ideals does not break off.}
\end{Example}

Although the similar result of Theorem 4.3 fails in the
Artinian case, we may prove the following 
Artinian version of Corollary 4.4, whose proof also holds in the
N\oe therian case.

\begin{Corollary} Let $N$ be a commutative algebra satisfying the identity $x^3=0$. If $N$ satisfies the descending
chain condition on ideals, then $N$ is finite-dimensional.
\end{Corollary}

{\it Proof. } When $N^2=0$, each subspace of $N$ is an ideal of $N$, and therefore $N$ must be finite-dimensional by the Artinian
hypothesis. Now, if $N^2\neq 0$, then the former case
shows that $N/N^2$ is finite-dimensional. Finally, we apply
\cite[Lemma 3.2]{Zitan} stating that any commutative algebra $N$
satisfying the identity $x^3=0$ and such that $N/N^2$ is
finite-dimensional must be finite-dimensional. \ep\\

At present, let us return to Bernstein algebras with locally nilpotent barideals. Applying Theorem
3.1 together with Theorem 4.3, we deduce immediately the next consequence:

\begin{Corollary} Let $A$ be a Bernstein algebra with locally nilpotent barideal $N$. If $A$ is N\oe therian, then $A$
is finite-dimensional.
\end{Corollary}

As attempted, we will give in the sequel a counter-example to the
Artinian version of Corollary 4.8. We emphasize that the locally
nilpotent algebras $N$ treated in Examples 4.5 and 4.6 cannot be
embedded in a Bernstein algebra $A$, since in either cases the identity
$\left(x^2\right)^2=0$ is not valid in $N$. For this reason, we
make appeal to the following example taken from
\cite[Example 3.10]{BZ}:

\begin{Example} {\rm Let $A$ be the Bernstein algebra with denumerable basis $\{e, v_1, u_1,
u_2, u_3, \dots\}$ and nonzero products
$$e^2=e,\ eu_i=\frac 12 u_i\ (i\geq 1),\ u_iv_1=u_{i-1} \ (i\geq
2)$$ The weight function $\omega: A \longrightarrow K$ is defined
by $\omega(e)=1,\ \omega(v_1)=\omega(u_i)=0 \ (i\geq 1)$ and the
Peirce components are $U=<u_1, u_2, \ldots> $ and $V=<v_1>$. We
know from \cite[Example 3.10]{BZ} that $A$ is Artinian. Actually,
it is clear that the barideal $N=U\oplus V$ is locally nilpotent.
We point out that the Bernstein algebra $A$ is not  N\oe therian, and moreover, it is neither
 nuclear nor Jordan.}
\end{Example}

\section{On the nilpotence}

\hspace*{0,5cm} In this final section we  shall proceed to revisit  some results of Micali and
Ouattara \cite{Mi}  in the aim to improve and generalize them to the N\oe therian and Artinian situations.  \\
First, we start with the following result which was proved in \cite[Lemmas
4.3 and 4.4]{Mi} when the Bernstein algebra $A$ was assumed to be finitely generated.

\begin{Lemma} Let $A=Ke \oplus U \oplus V$ be a Bernstein algebra which is N\oe therian or Artinian, and let $N=U\oplus V$ be its barideal. Let $I$ be a subspace of $A$. \\
	{\rm (i)} If $NI=I$, then $I\subseteq ann_U(U)$ and $I$ is an ideal of $A$.\\
	{\rm (ii)} $NI=I$ if and only if $VI=I$.
\end{Lemma}

{\it Proof. }  In view of \cite[Proposition 3.1]{BZ}, since $A$ is N\oe therian or Artinian, the Bernstein-Jordan algebra $A/ann_U(U)$ is finite-dimensional. Hence, its barideal $N/ann_U(U)$ is
nilpotent, so $N^k\subseteq ann_U(U)$ for some integer $k$. The remainder of the proof follows as in \cite[Lemmas 4.3 and 4.4]{Mi}. More precisely:\\
(i) Evidently, $I=NI\subseteq N$. Now, $I=NI=N(NI)=\dots \subseteq N^k$, yielding $I\subseteq ann_U(U)$. Furthermore, since $I\subseteq ann_U(U)\subseteq U$, we have $eI=I$, and by the condition $NI=I$, we deduce that $I$ is an ideal of $A$.\\
(ii) If $NI=I$, then the above assertion gives $I\subseteq ann_U(U)$. It follows that $I=NI=(U\oplus V)I=VI$, because $UI=0$. \\
Conversely, if $VI=I$, then $I=V(V(\dots V(VI)\dots))\subseteq N^k\subseteq ann_U(U)$. Therefore, $UI=0$ and so $NI=(U\oplus V)I=VI=I$. \ep\\

	\hspace*{0,5cm} Actually, after proving Lemma 5.1, the result of \cite[Th\'eor\`eme 4.7]{Mi} can
	be ameliorated by deleting the superfluous hypothesis that $A$ be
	finitely generated. Namely:
	
	\begin{Theorem} Let $A=Ke \oplus U \oplus V$ be an Artinian  Bernstein algebra. Then the following conditions are equivalent:\\
		{\rm (i)} The ideal $N=U\oplus V$ is nilpotent;\\
		{\rm (ii)} The associative algebra ${\cal M}_*^N(V)$ is nilpotent;\\
		{\rm (iii)} $I=0$ is the unique subspace of $A$ satisfying $VI=I$.
	\end{Theorem}
	
	{\it Proof. }  As in the proof of \cite[Th\'eor\`eme 4.7]{Mi}, the implications  ${\rm (i)} \Rightarrow {\rm (ii)} \Rightarrow {\rm (iii)}$ are always true even if $A$ is not Artinian. And the implication  (iii) $\Rightarrow$ (i) follows the same path as in \cite[Th\'eor\`eme 4.7]{Mi}, by applying our Lemma 5.1 instead of
	\cite[Lemma 4.3]{Mi}. In details:\\
	${\rm (i)}\Rightarrow {\rm (ii)}$: Since $N$ is nilpotent, then also is the multiplication ideal ${\cal M}_*(N)$ (see \cite[Chapitre II, Theorem 2.3]{Schafer}. In particular, the subalgebra ${\cal M}_*^N(V)$ is nilpotent.\\
	${\rm (ii)}\Rightarrow {\rm (iii)}$: Let $I$ be a subspace of $A$ with $VI=I$, so that $I\subseteq N$ by Lemma 5.1. Since $I=V(V(\dots V(VI)\dots))$ and ${\cal M}_*^N(V)$ is nilpotent, then $I=0$.\\
	${\rm (iii)}\Rightarrow {\rm (i)}$: Since $A$ is Artinian, the descending chain of ideals $(N^i)_{i\geq 1}$ of $A$ must stabilize. Hence, $N^r=N^{r+1}$ for some integer $r\geq 1$, or equivalently $N^r=NN^r$. It follows from Lemma 5.1(ii) that $N^r=VN^r$, implying $N^r=0$  by hypothesis. \ep
	
	\begin{Note}{\rm We point out that the implication ${\rm (ii)}\Rightarrow {\rm (i)}$ of
			Theorem 5.2 is already true when the Bernstein algebra $A$ is finitely
			generated \cite[Lemma 4.1]{Mi}, so it holds automatically in the N\oe therian
			case, since N\oe therian Bernstein algebras are finitely generated \cite[Corollary 3.6]{BZ}. Nevertheless, the implication ${\rm (iii)}
			\Rightarrow {\rm (i)}$ fails when $A$ is not Artinian, even if it is finitely generated or
			N\oe therian, as the following example will illustrate.}
	\end{Note} 
	
	\begin{Example} {\rm Let $A $ be the  infinite-dimensional Bernstein algebra considered in \cite[Example 3.11]{BZ}, with
			basis  $\{e, v_2, u_1, u_2, u_3, \ldots\} $ and nonzero
			products
			$$e^2=e,\; eu_i=\frac 12 u_i \; (i\geq 1), \;\; u_iv_2=u_{i+1} \; (i\geq 1),
			\; (i\geq 2). $$ Then $A=Ke\oplus U \oplus V $, where $U=<u_1, u_2,
			\ldots> $ and $V=<v_2> $.\\
			Let $I$ be a subspace of $A$ such that $VI=I$. Assume that $I\neq 0$ and choose an element
			$a=(\alpha_1u_1+\dots+\alpha_{p}u_p)+\alpha v_2 \in I$, with $p$ minimal such that $\alpha_p\neq 0$. Since $a\in VI$, there exists
			$b=(\beta_1u_1+\dots+\beta_{q}u_q)+\beta v_2\in I$ with $\beta_q\neq0$ and $a=v_2b=\beta_1u_2+\dots+\beta_qu_{q+1}$. Then the
			contradiction $q+1=p<q$ yields $I=0$. However, the ideal $N=U\oplus V$ is not nilpotent. In fact, this Bernstein algebra $A$ is finitely generated, N\oe therian but not Artinian  \cite[Example 3.11]{BZ}.}
		\end{Example}
	
	 \quad We close our paper by making the following  comment. The Grishkov conjecture \cite{Gr} asserts that if $A=Ke\oplus U\oplus V$ is a finitely generated Bernstein algebra that is nuclear, then the barideal $N=U\oplus V$ is nilpotent. This question has been shown affirmatively by  Peresi \cite{Pe} and  Krapivin \cite{Kra}. Nevertheless, the  proof of this result presented in \cite[Th\'eor\`eme 4.10]{Mi}  is not correct,
because it relies on \cite[Th\'eor\`eme 4.7]{Mi} which requires the additional assumption that $A$ be Artinian.
$$$$

	{\bf \large Acknowledgment}: 	
	The author wish to thank Professor Nadia Boudi for a number of useful  discussions  and for reading an earlier draft of this paper.

\end{document}